\documentclass[11pt,reqno]{amsart}
\usepackage{epsfig}

\usepackage{amssymb,amsmath}

\newcommand{\be}{\begin{eqnarray}}
\newcommand{\ee}{\end{eqnarray}}

\newcommand{\const}{\mbox{\rm const}}

\newcommand{\es}{\emptyset}

\newcommand{\eps}{{\mbox{$\epsilon$}}}

\newcommand{\R}{{\mathbb R}}
\newcommand{\Q}{{\mathbb Q}}
\newcommand{\Z}{{\mathbb Z}}

\newcommand{\Nat}{{\mathbb N}}

\newcommand{\Ak}{{\mathcal A}}

\newcommand{\Lk}{{\mathcal L}}

\newcommand{\Lam}{{\Lambda}}

\newcommand{\sig}{\sigma}

\newcommand{\Sk}{{\mathcal S}}
\newcommand{\Tk}{{\mathcal T}}
\newcommand{\Xt}{X_\Tk}
\newcommand{\bt}{{\bf t}}
\newcommand{\piga}{\pi_\gamma}
\newcommand{\pigas}{\pi_\gamma^\perp}
\newcommand{\old}[1]{}

\def\th{\theta}
\def\gam{\gamma}
\def\Gam{\Gamma}
\def\Phti{\widetilde{\Phi}}
\def\Gati{\widetilde{\Gamma}}
\def\La{\Lam}
\def\Ga{\Gam}

\newtheorem{theorem}{Theorem}[section]
\newtheorem{lemma}[theorem]{Lemma}
\newtheorem{cor}[theorem]{Corollary}
\newtheorem{prop}[theorem]{Proposition}

\theoremstyle{definition}
\newtheorem{defi}[theorem]{Definition}
\newtheorem{example}[theorem]{Example}

\theoremstyle{remark}

\numberwithin{equation}{section}

\begin{document}

\title[Topological mixing for substitution systems]
{Topological mixing 
for substitutions on two letters}

\author{Richard Kenyon}
\address{Richard Kenyon, Department of Mathematics, University of British Columbia,
Vancouver, BC, V6T 2Z1 Canada}
\email{kenyon@math.ubc.ca}

\author{Lorenzo Sadun}
\address{Lorenzo Sadun, Department of Mathematics, University of Texas,
Austin, TX 78712}
\email{sadun@math.utexas.edu}
\thanks{Research of Sadun was supported in part by NSF grant \#DMS 0305505}

\author{Boris Solomyak}

\thanks{
Research of Solomyak was supported in part by NSF grant
\#DMS 0099814} 

\address{Boris Solomyak, Box 354350, Department of Mathematics,
University of Washington, Seattle, WA 98195}
\email{solomyak@math.washington.edu}


\begin{abstract} We investigate topological mixing for $\Z$ and $\R$ actions
associated with primitive substitutions on two letters. The characterization
is complete if the second eigenvalue $\theta_2$ of the substitution matrix
satisfies $|\theta_2|\ne 1$. If $|\theta_2|<1$, then (as is well-known) the
substitution system is not topologically weak mixing, so it is not
topologically mixing. We prove that if $|\theta_2|> 1$, then
topological mixing is equivalent to topological  weak mixing, which
has an explicit arithmetic characterization. The case $|\theta_2|=1$ is
more delicate, and we only obtain some partial results.
\end{abstract}

\maketitle

\thispagestyle{empty}
\section{Introduction and statement of results}

Let $X$ be a compact metric space and let $G=\Z^d$ or $\R^d$
act continuously on $X$. Let $|g|$ be the
distance from $g$ to $0$ in some translation-invariant metric.

The dynamical system $(X,G)$ is said to be {\em topologically mixing} if
for any two nonempty open sets $U, V \subset X$, there exists $R>0$ such that
$$
U \cap T_g(V) \ne \es,\ \ \ \mbox{for all}\ g\in G,\ |g| \ge R.
$$
The dynamical system is {\em topologically weak mixing} if it has no
non-constant continuous eigenfunctions. It is easy to see that 
topological mixing implies topological weak mixing.

\old{It is well-known that topological mixing implies that the dynamical system has
no nonconstant continuous eigenfunctions 
(the latter property is called topological weak
mixing).}

For a symbolic dynamical system, 
topological mixing is equivalent to the property that
for any two allowed blocks $W_1$ and $W_2$ 
there exists $N\in \Nat$ such that for any
$n\ge N$ there is an allowed block $W_1BW_2$, with $|B|=n$.

\medskip

Let $\Tk$ be a tiling of $\R_+$ with a finite set of interval
prototiles, and let $\Xt$ be the associated tiling space.
That is, $\Xt$ is the set of all tilings $\Sk$ of $\R$
such that every patch of $\Sk$ is the translate of a patch in $\Tk$.
(A {\em patch} is a tiling of a finite interval, and two patches are 
said to be {\em equivalent} if each is a translate of the other).  
The {\em tiling dynamical system} is the $\R$-action $(\Xt,T_g)$
where $T_g(\Sk)=\Sk-g$.
The topology on $\Xt$ has basis given by ``cylinder sets'' indexed by 
a patch $P$ and a radius $\eps$:
$$
X_{P,\eps} = \{\Sk \in \Xt:\ \exists y \in (-\eps,\eps):  \ P-y \subset
\Sk\}.
$$

It follows that $(\Xt,T_g)$ is topologically mixing if and only if for
any two patches $P_1, P_2 \subset \Tk$, and for any $\eps>0$, there
exists $R>0$ such that for all $g\in \R$, with $|g|\ge R$, there exist
$\Sk\in \Xt$ and $y \in (-\eps,\eps)$, such that
$$
P_1 \subset \Sk\ \ \ \mbox{and}\ \ \ P_2 - g - y \subset \Sk.
$$

We say that a set $Y\subset \R$ is {\em eventually dense} in $\R$ if
for any $\eps>0$ there exists $R>0$ such that
the $\eps$-neighborhood of $Y$ covers $\R \setminus
(-R,R)$. An alternative way to state topological mixing for a tiling system
is to say that for any allowed patches $P_1$ and $P_2$,
the set of translation vectors between locations (say,
left endpoints) of patches equivalent to $P_1$ and patches equivalent to $P_2$
is eventually dense in $\R$. 

We will also say that $Y \subset \R_+$ is eventually dense in $\R_+$ if
$Y \cup -Y$ is eventually dense in $\R$.

\medskip

Now we briefly recall the definition of substitutions and associated
dynamical systems; see \cite{queffelec,Pyth} for more details.
Let $\zeta$ be a substitution on a finite alphabet $\Ak = 
\{0,1,\ldots,m-1\}$, with $m\ge 2$. Recall that $\zeta$ is a mapping from
$\Ak$ to $\Ak^*$, the set of nonempty words in the alphabet $\Ak$.
The substitution $\zeta$ is extended to maps (also denoted by $\zeta$)
$\Ak^*\to\Ak^*$ and $\Ak^{\Nat} \to\Ak^{\Nat}$ by concatenation.
The matrix $M=M_\zeta$ 
associated with the substitution $\zeta$ is the $m\times m$ matrix
defined by
$$
M = (m_{i,j})_{d\times d},\ \ \ \mbox{where}\ \ \
m_{i,j} = \ell_i(\zeta(j))
$$
where $\ell_i(W)$ is the number of occurrences of letter $i$ in word $W$.
We will always assume that the
substitution is {\em primitive}, that is, there exists $k$ such that
all entries of $M^k$ are strictly positive; equivalently,
for every $i,j \in \Ak,$ the symbol
$j$ occurs in $\zeta^k(i)$. The {\em substitution space} $X_\zeta$ is defined
as the set of all two-sided infinite sequences in the alphabet $\Ak$ whose
every word occurs in $\zeta^k(i)$ for $k$ sufficiently large. Then
$X_\zeta$ is a closed shift-invariant subset of $\Ak^\Z$;
the $\Z$-action associated with the substitution is $(X_\zeta,\sig)$
where $\sig$ is the left shift. We may assume, without loss of generality,
that $\zeta(0)$ starts with $0$ (if not, replace $\zeta$ with $\zeta^k$
for appropriate $k$ and rename the symbols). Then we get a one-sided
fixed point of the substitution map $u = \zeta(u) = \lim_{n\to \infty}
\zeta^n(0)$, sometimes called the {\em substitution sequence}. 
Primitive substitution $\Z$-actions are uniquely ergodic and minimal, 
see \cite{queffelec}, 
which implies that every
allowed word in the substitution space occurs in $u$. 
We denote by $\Lk(X_\zeta)$ the
{\em language of the substitution}, that is the collection of all allowed
words (or equivalently, all subwords of $u$).
In addition, we will always assume that 
$\zeta$ is {\em aperiodic}, that is, the substitution sequence
$u$ is  not periodic, which is equivalent to $X_\zeta$ being infinite.

\medskip

Next we recall the definition of the tiling dynamical system
associated with the substitution $\zeta$; see \cite{Soltil,CS} for
more details.  Let $\bt = (t_i)_{i\in \Ak}$ be a strictly positive row
vector.  We associate to $u$ a tiling $\Tk$ of $\R_+$ whose prototiles
are intervals $\tau_i$ of length $t_i$, for $i=0,\ldots,m-1$, in such
a way that $0$ is the left endpoint of the tile $\tau_{u_1}$, followed
by the copy of $\tau_{u_2}$, etc.  {\em A priori}, some of the tiles
may be congruent; then we distinguish them by ``labels'' from the
alphabet $\Ak$.  The tiling space $\Xt$ is defined, as above, as the
set of tilings $\Sk$ of the line $\R$ such that every patch of $\Sk$
is a translate of a $\Tk$-patch, and the group $\R$ acts by
translation.  Observe that this $\R$-action is topologically conjugate
to the suspension flow over the $\Z$-action $(X_{\zeta},\sigma)$, with
the height function equal to $t_i$ on the cylinder corresponding to
the symbol $i$.  This system is also minimal and uniquely ergodic.

\medskip

Mixing properties of general
uniquely ergodic systems have been much investigated,
and we do not survey this literature here. Many early references can be found
in the paper by Petersen and Shapiro \cite{PeSh}. For such systems,
measure-theoretic strong mixing implies topological strong mixing, which
implies topological weak mixing, and none of the implications can be 
reversed, see \cite{PeSh}. Mixing properties of primitive
substitution $\Z$-actions were studied by Dekking and Keane \cite{DK},
who proved that they are never strongly mixing, but may be topologically
mixing. More recently, topological mixing for substitutions on two symbols
was investigated by A. Livshits \cite{liv1,liv2,liv3}, but a characterization
of such systems was still lacking. In another, though related, direction,
Host \cite{Host} proved that for substitution $\Z$-actions
topological weak mixing is equivalent to measure-theoretic weak mixing.
For substitution $\R$-actions this is essentially proved in
\cite[Thm 2.3]{CS}. (Alternatively, one can argue as in
\cite[Thm 4.3]{Soltil}, using \cite[Lem 2.1]{CS} instead of
\cite[Lem 4.2]{Soltil}.)
Thus, when writing ``weak mixing'' we will not specify whether it is
in the topological or measure-theoretic category. 
 
\medskip

We begin with a standard elementary proposition which contains
general necessary conditions for topological mixing of substitution systems.

\begin{prop} \label{prop-nec}
Let $\zeta$ be a primitive aperiodic substitution on the alphabet $\Ak$.

{\bf (i)} If the $\Z$-action $(X_\zeta,\sig)$ is topologically mixing, then
\be \label{reprime}
GCD\{|\zeta^n(i)|:\ i \in \Ak\} = 1,\ \ \ \mbox{for all}\ n\ge 1.
\ee

{\bf (ii)} Let $(t_i)_{i\in \Ak}$ be a vector of tile lengths.
If the $\R$-action $(\Xt,T_g)$ is topologically mixing, then
\be \label{irrat}
\mbox{There exist } i,j \in \Ak \mbox{ such that }
t_i/t_j \mbox{ is irrational.}
\ee
\end{prop}

{\em Proof.}
(i) Suppose that for some $m\in \Nat$, the numbers
$|\zeta^m(i)|$,\ $i\in \Ak$, have a common factor $p$.
Primitive aperiodic substitutions are known to be bilaterally recognizable
\cite{mosse}. This implies that if an allowed word $W$ is sufficiently long,
then it is uniquely determined where the subwords $\zeta^m(i)$ occur within
$W$ (after some ``buffer'' of uniformly bounded length is 
removed from the beginning and the end) in any occurrence of $W$ in the
substitution sequence $u$. Then it follows that the distance between
any two occurrences of $W$ is divisible by $p$, contradicting topological
mixing.

\medskip

(ii) Clearly, the left endpoints of all tiles of the tiling $\Tk$ belong to 
the $\Z$-module generated by the $t_i$ where $\ i\in \Ak$. 
If the lengths are all rationally related, then
this $\Z$-module is a discrete subset of $\R$, which contradicts topological
mixing for the $\R$-action.
\qed

\medskip

 From now on we consider primitive aperiodic substitutions on 2 symbols only.
We will use the basic facts of Perron-Frobenius theory, see, e.g.\ 
\cite{Sen}.
Let $\theta_1$ be the Perron-Frobenius eigenvalue of the substitution matrix
$M$ and let $\theta_2$ be the second eigenvalue.
Our main result is

\begin{theorem} \label{th-main}
Suppose that $|\theta_2| > 1$. Then

{\bf (i)} the $\Z$-action $(X_\zeta,\sig)$ is topologically mixing 
if and only if (\ref{reprime}) holds;

{\bf (ii)} the $\R$-action $(\Xt,T_g)$ is topologically mixing if and only if
(\ref{irrat}) holds (that is, if $t_1/t_0 \not\in \Q$).
\end{theorem}

\noindent{\bf Remarks.} 

1. If $|\theta_2|<1$, then both the $\Z$-action and the
$\R$-action associated with the substitution have nontrivial 
continuous eigenfunctions \cite{Host,CS},
so they are not topologically mixing. The case $|\theta_2|=1$ is more subtle.
Dekking and Keane \cite{DK} considered the following two substitutions:
$$
\zeta_1(0) = 001,\ \ \zeta_1(1) = 11100;\ \ \ \ \ 
\zeta_2(0) = 001,\ \ \zeta_2(1) = 11001.
$$
They have the same substitution matrix with eigenvalues $4$ and $1$.
In \cite{DK} it is proved that the $\Z$-action associated with $\zeta_1$
is topologically mixing, whereas the one associated with $\zeta_2$ is not
(the latter actually goes back to \cite{PeSh}). We show (Theorem \ref{thm2}(ii) and
section \ref{examplesection}) that the same
is true for the corresponding $\R$-actions 
(with irrational ratio of tile lengths).

2. Partial results in the direction of conclusion (i) were obtained by 
A. Livshits \cite{liv1,liv2,liv3}, and in fact some of our methods are 
similar to his. He conjectured that (i) holds.

3. As far as we are aware, topological mixing for 
substitution $\R$-actions has not been
considered before. A special choice of the tile lengths is 
the one arising from the
Perron-Frobenius eigenvector of the substitution matrix. 
Then we get a (geometrically) self-similar
tiling of the half-line. 

4. Given a substitution, it is straighforward to check condition
(\ref{reprime}). Let $r = (1,\ldots,1)$.  Then $r M^n =
(|\zeta^n(0)|,\ldots,|\zeta^n(m-1)|)$.  The question is whether, for
any prime $p$, this vector is eventually zero mod ${p}$.  If $p$ does
not divide the determinant of $M$, then $M$ is invertible mod ${p}$,
and $r M^n$ can never equal zero mod ${p}$. Thus we need only consider
primes that divide the determinant of $M$. For each of these, the
sequence of vectors $\{ rM^n \pmod{p} \}$ takes on at most $p^{m}$ values, 
hence starts repeating after at
most $p^{m}$ terms, so we need only examine those first $p^m$ terms.

5.  Results similar to Theorem \ref{th-main} are known for weak mixing 
\cite{Host,CS}. The difference is that
$|\theta_2| \ge 1$ (as opposed to $|\theta_2| > 1$), 
together with (\ref{reprime}) or (\ref{irrat}), implies weak mixing. Thus

\begin{cor} \label{cor} Consider a tiling
dynamical system arising from a primitive substitution on two 
symbols. If $|\theta_2| \ne 1$, then the $\R$ action is topologically mixing 
if and only if it is weak mixing. 
\end{cor}

\medskip Although a complete description of substitutions with
$|\theta_2|=1$ remains open, there is something
that we can say. The following theorem
constitutes one of the two main steps in the proof of Theorem~\ref{th-main},
and it deals with an arbitrary primitive substitution on two symbols.
Consider all words $W$ in $\Lk(X_\zeta)$ of length $|W|=n$. Let
$$ a(n) = \min \ell_0(W), \qquad \qquad b(n) = \max \ell_0(W)$$
be the minimum and maximum
number of 0's in allowed words of length $n$. We will need the condition
\be \label{excess}
b(n) - a(n) \to \infty,\ \ \ \mbox{as}\ n\to \infty,
\ee
which is closely related to the ``growth of excess" condition that appeared in
\cite{DK} and \cite{liv1,liv2,liv3}.

\begin{theorem} \label{thm2}

{\bf (i)} If (\ref{reprime}) is satisfied, then
the $\Z$-action corresponding to the
substitution is topologically mixing if and only if
(\ref{excess}) holds.

{\bf (ii)} If (\ref{irrat}) is satisfied, then
the $\R$-action corresponding to the
substitution is topologically mixing if and only if 
(\ref{excess}) holds.
\end{theorem}

\medskip

The next proposition is the second ingredient of the proof of 
Theorem~\ref{th-main}.

\begin{prop} \label{prop1}
Suppose that the substitution satisfies the condition $|\theta_2|>1$.
Then there exists a constant $c_1>0$ such that 
$$
b(n) - a(n) \ge c_1 n^\alpha,\ \ \ \mbox{where}\ 
\alpha = \log|\theta_2|/\log\theta_1 \in (0,1).
$$
\end{prop}

It is clear that Theorem~\ref{th-main}
will follow once we prove Theorem \ref{thm2} and
Proposition \ref{prop1}.

The difference between topological mixing and weak mixing is related
to the difference between lim-sup and lim-inf.  Consider the
quantities $\sup_{n<N}(b(n)-a(n))$ and $\inf_{n>N} (b(n) - a(n))$. If
$|\theta_2|>1$, then both of these quantities grow as $N^\alpha$
and the system is both topologically mixing and weak mixing.  If
$|\theta_2|<1$, then both of these quantities are bounded and the
system is neither topologically mixing nor weak mixing.  If $|\theta_2|=1$,
then, as can be shown,
the sup grows as $\log(N)$ (implying weak mixing) but the inf may
or may not grow.  

\section{Preliminaries}

For a word $W$ consisting of 0's and 1's, recall that
$\ell_i(W)$ denotes the number of letters
$i$ in $W$. The column vector $\ell(W) = 
(\ell_0(W), \ell_1(W))^T$ (where ${}^T$ denotes transpose) 
is called the {\em population vector} of the word $W$.
By the definition of the substitution matrix, for any word $V$ in the
alphabet $\Ak$,
\be \label{matr}
\ell(\zeta(V)) = M \, \ell(V).
\ee
Note that the length of a word $V$ is $|V|=(1,1)\cdot \ell(V)$. 
In a tiling, the length of a patch corresponding to the word $V$ (also called
the {\em tiling length} of $V$ and denoted $|V|_\Tk$) is 
$(t_0, t_1)\cdot \ell(V)$.  Recall that $\Lk(X_\zeta)$ denotes
the {\em language} of the substitution subshift,
that is, the set of all words that occur in the subshift.
Consider
$$
\Phi(X_\zeta) := \{\ell(W):\ W\in \Lk(X_\zeta)\}.
$$
Further, consider the set $\Gamma(X_\zeta)$ of all the points $\ell(u[1,j])$
for 
$j \ge 0$ (for $j=0$ we just get the origin). 
Observe that 
$$
\Phi(X_\zeta) = (\Gamma(X_\zeta)-\Gamma(X_\zeta))\cap \Z_+^2.
$$

When the substitution is
fixed, we drop $X_\zeta$ from the notation and write just $\Lk, \Phi, \Gamma$.
Connecting the consecutive points of $\Gam$ we obtain a polygonal curve
$\Gati$, which gives a nice visual representation of the substitution sequence
$u$. This curve starts at
the origin and goes into the 1st quadrant, with edges going up or to the
right along
the standard grid. Clearly $(1, 0)^T
\in \Gamma$ since $u$ starts with 0.

\begin{lemma} \label{lem-elem1}

{\bf (i)} $\Phi = \bigcup_{n\ge 1} \{(i, n-i)^T
:\ a(n) \le i \le b(n)\}$.

{\bf (ii)} $0 \le b(n+1) - b(n) \le 1$;\ \
$0 \le a(n+1) - a(n) \le 1$.

{\bf (iii)} If $(i, j)^T,
(i',j')^T \in \Phi$ and 
$i \le i', j\ge j'$, then $(k, m)^T \in
\Phi$ whenever $i \le k \le i'$ and $j' \le m \le j$.
\end{lemma}

{\em Proof.} (i) Consider the population vectors of all words of length $n$.
They include $(a(n), n-a(n))^T$ and $(b(n), n-b(n))^T$ and hence all 
intermediate ones since $|\ell_0(u[i, i+n-1]) - \ell_0(u[i+1,i+n])| \le 1$.

(ii) is obvious from the definition.

(iii) follows easily from (i) and (ii). The details are left as an exercise.
\qed

\medskip

By Lemma~\ref{lem-elem1}(i),
the ``upper envelope'' of $\Phi$ is the set $\{(a(n), n-a(n)):\ n\ge 1\}\cup
(0,0)$, and the ``lower  envelope'' of $\Phi$ is the set $\{(b(n), n-b(n)):\
 n\ge 1\}\cup (0,0)$. We connect the consecutive points of these ``envelopes''
to obtain two polygonal curves and consider the set $\Phti$ between
them. It is a kind of ``strip'' with polygonal edges; by definition,
$\Phi = \Phti\cap \Z^2_+$.

\medskip

\begin{example} Let $\zeta(0)= 011,\ \zeta(1)=0$. The matrix of the 
substitution is $M= \left[ \begin{array}{ll} 1 & 1 \\ 2 & 0 
\end{array} \right]$, with the eigenvalues $\theta_1 = 2$, $\theta_2 = -1$.
The fixed point of the substitution is 
$$
u = 0110001101101100011000110001101101100011011...
$$
This fixed point is represented by $\Gati$, shown as the thick line in
Figure 1. The region $\Phti$ is shaded. We will show in Section 
\ref{examplesection}
that for this substitution $\sup(b(n)-a(n)) = \infty$ but $\liminf
(b(n)-a(n)) \le 2$.
\end{example}

\begin{figure}[ht]
\epsfig{figure=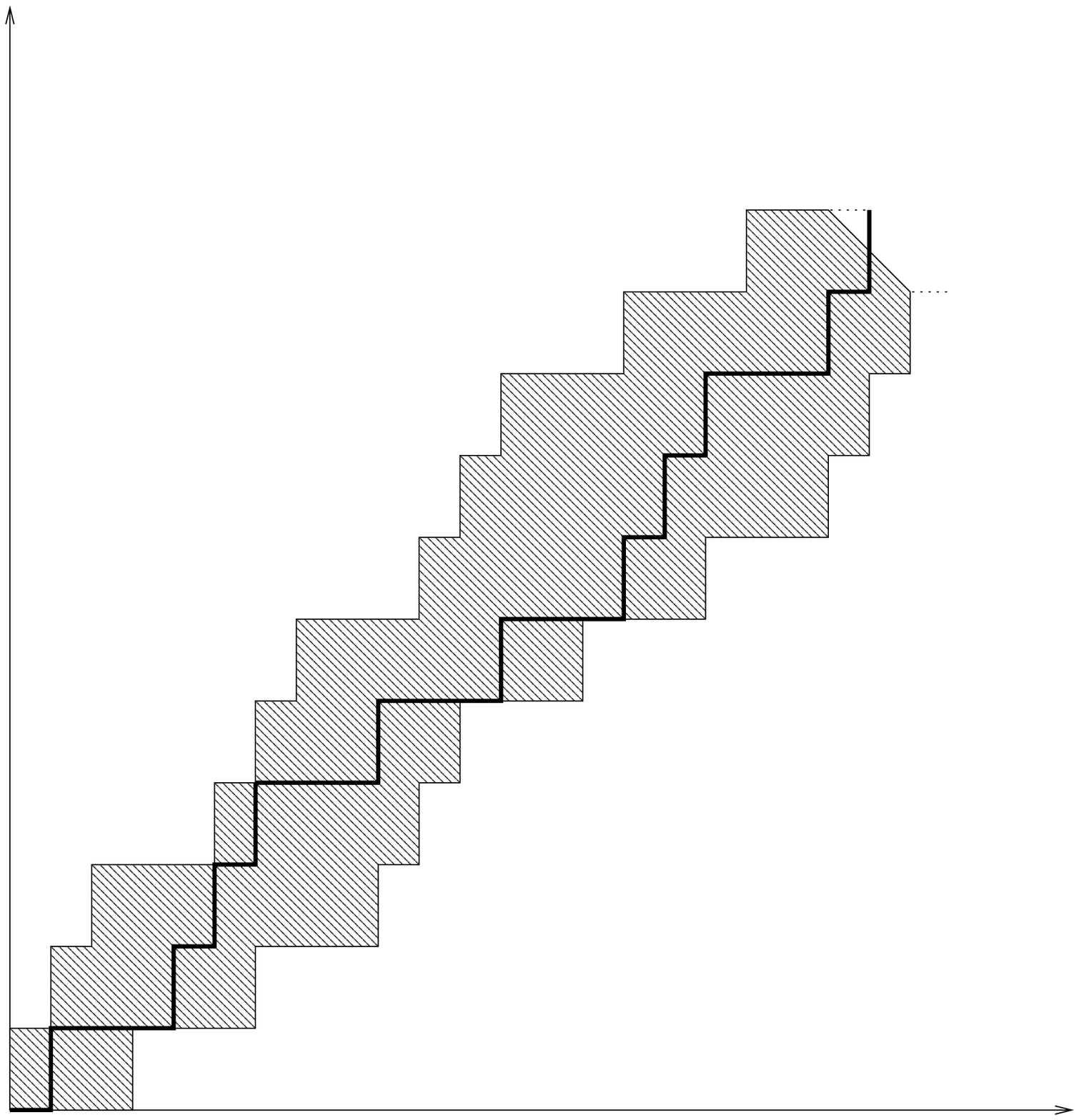,height=14cm}
\caption{The sets $\Gati$ and $\Phti$}
\end{figure}

\begin{defi} \label{def-width}
For any $\gam>0$, we define the width of $\Phti$
in the direction of
$(-\gam,1)$ at the level $r>0$ to be the length of the intersection of
$\Phti$ with the line $x+\gam y = r$. 
\end{defi}

For a tiling system with ${\bf t}=(1,\gamma)$, the width of $\Phti$
in the direction of $(-\gam,1)$ at the level $r$ governs the number of 
distinct population vectors whose tiling lengths are approximately $r$. 
In particular, for $\gamma$ irrational, 
the set of tiling lengths of allowed words 
is asymptotically dense if and only if the 
width of $\Phti$ goes to infinity as $r \to \infty$.  This geometric
observation is at the heart of the proofs of Propositions 3.3 and 3.5, 
below. 

Observe that $\sqrt{2}(b(n)-a(n))$ is the width of $\Phti$
in the direction of $(-1,1)$ at the level $n$ (since the points
$(a(n), n-a(n))^T$ and $(b(n), n-b(n))^T$ lie on the line $x+y=n$).  
The following lemma 
shows that the asymptotic behavior of the width for large levels does not
depend on the direction, and hence is determined by the large-$n$
behavior of $b(n)-a(n)$. 

\begin{lemma} \label{lem-elem22} For any $\gam_1,\gam_2>0$ there exist
constants $K_{\gam_1,\gam_2},C_{\gam_1,\gam_2}>1$ such that if 
$\Phti$ has width $L_1$ in the direction of
$(-\gam_1,1)^T$ at the level $r_1$, then there is 
$$
r_2 \in (K_{\gam_1,\gam_2}^{-1} r_1, K_{\gam_1,\gam_2} r_1) 
$$
such that the width $L_2$ of $\Phti$ in the direction of $(-\gam_2,1)^T$ 
at the level $r_2$ satisfies
$$
C^{-1}_{\gam_1,\gam_2} L_1 \le L_2 \le C_{\gam_1,\gam_2} L_1.
$$
\end{lemma}

{\em Proof.} 
This is a simple geometric fact. From
Lemma~\ref{lem-elem1}(iii), 
$\tilde\Phi$ contains a rectangle $R_1$ with sides parallel to the axes,
with aspect ratio $\gamma_1$ and whose diagonal is 
on the line $x+\gamma_1 y=r_1$ and has length $L_1$.
For a well-chosen $r_2$, there is a rectangle $R_2$
of aspect ratio $\gamma_2$ contained in this rectangle,
whose diagonal is on the line $x+\gamma_2y=r_2$, and whose diagonal length
is proportional to $L_1$.
Clearly $r_2\in(K^{-1}r_1,Kr_1)$
for some constant $K$ depending only on $\gamma_1$ and $\gamma_2$.
\qed

\medskip

Let $e_i$ be eigenvectors of $M$ corresponding to $\theta_i$, for
$i=1,2$. By the Perron-Frobenius theory, $e_1$ has components of the same
sign and $e_2$ has components of opposite signs. We choose
$e_1$ to be strictly positive and $e_2$ such that $(e_2)_x < 0$. 
Denote
$$\gamma:=(e_1)_y/(e_1)_x,\  \ \ \ 
\alpha = \log|\theta_2|/\log\theta_1.
$$

For $w = (w_x, w_y) \in \R^2$, we define the maps
\begin{equation}\label{projs}
\piga(w) = w_x + \gamma w_y,\ \ \ \
\pigas(w) = w_y - \gamma w_x.
\end{equation}
Up to an overall factor of $\sqrt{1+\gamma^2}$, these give the (signed) length
of the projections of $w$ onto lines parallel and perpendicular
to $e_1$, respectively.
Note that $\pigas(e_1)=0$ and $\pigas(e_2)>0$.


\section{Proof of Theorem \ref{thm2}}

In this section we prove a series of propositions.  The first three show
that (\ref{excess}) is a necessary condition for topological mixing.
The last three show that (\ref{excess}), together with (\ref{reprime})
or (\ref{irrat}), implies mixing. The next proposition follows
from \cite[Theorem 22]{Adam}, which contains precise bounds
on $\limsup (b(n)-a(n))$. We provide a proof for completeness.

\begin{prop} \label{prop-exc1}
Suppose that a primitive aperiodic
substitution $\zeta$ on two letters has the second eigenvalue
$\theta_2$ satisfying $|\theta_2| \ge 1$. 
Then $\limsup (b(n)-a(n)) = \infty$.
\end{prop}

{\em Proof. } 
For any $n$, $b(n) \ge n \gamma/(\gamma + 1) \ge a(n)$, since
the average density of 0's is $\gamma/(\gamma+1)$.  To show that 
$\limsup (b(n)-a(n)) = \infty$, it suffices to show that either
$b(n) - \frac{n \gamma}{\gamma+1}$ or $\frac{n \gamma}{\gamma+1} - a(n)$
is unbounded.  This is 
equivalent to showing that there are vectors $w \in \Gamma-\Gamma$ with
$|\pigas(w)|$ arbitrarily large.

For $|\theta_2|>1$ this is easy, since $\pigas(Mw) = \theta_2 \pigas(w)$
and since $\pigas(1,0)^T \ne 0$. 
All that remains are the cases 
$\theta_2 = \pm 1$.  By squaring the substitution, we can assume that
$\theta_2=1$.

We can choose $i\in\Ak$ so that $ii$ is allowed, since
otherwise the substitution is periodic.
Then, by primitivity,
there exists $p\in \Nat$ such that $\zeta^p(i) = V_1ii  V_2$ for
some words $V_1, V_2$. Note that at least one of $\ell(V_2), \ell(iV_2)$ is
not a multiple of $e_1$. Thus we can write
$$
\zeta^p(i) = W_1 i W_2
$$
where 
$\ell(W_2) = a_1e_1 + a_2 e_2$, with $a_2 \ne 0$. Then
$$
\zeta^{2p}(i) = \zeta^p(W_1) W_1 i W_2 \zeta^p(W_2),
$$
and iterating this we obtain that
$$
\forall \,m\ge 1,\ U_m:= W_2\zeta^p(W_2)\ldots \zeta^{mp}(W_2) \in \Lk(X_\zeta).
$$
It follows that
\begin{equation} \label{kar2}
\pigas(\ell(U_m)) = a_2 \pigas(e_2) m,
\end{equation}
which can be made arbitrarily large.  \qed

\begin{prop} \label{prop-gen}
Let $S$ be any subshift on 2 letters (not necessarily generated 
by a substitution).  Let $b(n)$ be the maximum number of 0's in words of 
length $n$, and let $a(n)$ be the minimum.  If $\sup(b(n)-a(n)) > \liminf 
(b(n) - a(n))$, then the $\Z$ action on $S$ is not topologically mixing.
\end{prop}

\medskip

{\em Proof.}
Let $V_n$  (resp. $W_n$) be a word of length $n$ that maximizes (resp. 
minimizes) the number of 0's.  Suppose that $m>n$ and that $b(m)-a(m)$ is 
strictly less than $b(n)-a(n)$. Then words of the form $V_n U W_n$, where $U$ 
is a spacer of size $m-n$, can never occur, as $V_n U$ and $U W_m$ differ in 
population by more than $b(m)-a(m)$.  That is, the cylinder set based on 
$V_n$, and the translate by $m$ of the cylinder set based on $W_n$, do not 
intersect.

Now pick $n$ such that $b(n) - a(n)$ is larger than the $\liminf$.
Since there are arbitrarily large values of $m$ for which 
$b(n)-a(n) > b(m)-a(m)$, the cylinders on $V_n$ and $W_n$ do not mix. \qed

\medskip

\begin{sloppypar}
\begin{prop} \label{nec-tiling}  
  Let $\Xt$ be a tiling space based on a substitution on two letters.
  If $\liminf \left (b(n)-a(n)\right ) < \infty$, then the $\R$ action
  on $\Xt$ is not topologically mixing.
\end{prop}
\end{sloppypar}

{\em Proof.}
We can assume without loss of generality that the tile lengths are $t_0=1,
\ t_1=\beta>0$. Then $\pi_\beta(\Phi)$ is the set of tiling lengths of
words in $\Lk(X_\zeta)$ 
(here $\pi_\beta$ is defined by $\pi_\beta(w_x,w_y)=w_x+\beta w_y$
as in \eqref{projs}).
It is enough to show that 
$\pi_{\beta}(\Phi)$ is not eventually dense. 
Since $\liminf(b(n)-a(n)) < \infty$, 
the width of the set $\Phti$ in the direction $(-1,1)$ is bounded
at a sequence of levels $k_n \to \infty$. Then the width of $\Phti$ in
the direction $(-\beta, 1)$ is bounded at a sequence of levels
$k_n' \to \infty$, by Lemma~\ref{lem-elem22}. Since $\Phi$ is the set of
lattice points in $\Phti$, it follows that $\pi_\beta(\Phi)$
is not eventually dense. \qed

\medskip

Propositions (\ref{prop-exc1}, \ref{prop-gen}, \ref{nec-tiling}) imply half of 
Theorem \ref{thm2}, namely that (\ref{excess}) is a necessary condition
for topological mixing.  The remainder of this section is to prove
sufficiency. 

Denote
\be \label{Phipr}
\Phi' = \{\ell(0W):\ 0W1 \in \Lk\}.
\ee

\begin{lemma} \label{lem-elem2} 
$$
\Phi \setminus \{(a(n), n-a(n)):\ n \ge 1\} \subset \Phi'.
$$
\end{lemma}

{\em Proof.} 
Fix $n$ and consider $f_n(j) = \ell_0(u[j, j+n-1])$.
We have $|f_n(j+1)-f_n(j)| \le 1$, with $f_n(j+1)-f_n(j)=-1$ if and only if
$u_j=0$ and $u_{j+n}=1$. Since $u$ is uniformly recurrent (the substitution
being primitive), $a(n) = \min_j f_n(j)$ and $b(n) = \max_j f_n(j)$ are
achieved infinitely often. It follows that
\be \label{imp1}
\forall k\in [a(n) +1, b(n)] \cap \Z,\ \exists j\ge 1:\ \ u_j=0, u_{j+n} = 1,
\ \ \ell_0(u[j,j+n-1])=k.
\ee
This implies the desired statement, in view of Lemma \ref{lem-elem1}. \qed

\medskip

\begin{prop} \label{suff-1}
Let $\Xt$ be a tiling space based on a substitution on two letters.
If $\liminf (b(n)-a(n)) = \infty$ and the ratio of tile lengths is 
irrational, then the $\R$ action on $T$ is topologically mixing.
\end{prop}

{\em Proof.}
By the definition of our tiling space, all patches are determined by allowed
words of the substitution space. Recall that the tiling length of a word
$V$ is 
$$
|V|_{\Tk} = \ell_0(V) t_0 + \ell_1(V) t_1 = \bt \cdot \ell(V),
$$
where $\bt = (t_0,t_1)$ is the vector of prototile lengths. 
To prove topological mixing, we
need to show that for any allowed words $W_1,W_2$, the set 
$$
\Psi(W_1,W_2) := \{|W_1V|_{\Tk}:\ W_1 V W_2 \in \Lk\}
$$
is eventually dense in $\R_+$.
Since $\zeta$ is a primitive substitution,
there exists $m$ such that $\zeta^m(0)$ contains $W_1$ and $\zeta^m(1)$ contains
$W_2$. Thus, it is enough to prove the claim for $W_1 = \zeta^m(0)$ and
$W_2 = \zeta^m(1)$ for an arbitrary $m\in \Nat$.
Clearly,
$$
\Psi(\zeta^m(0),\zeta^m(1)) \supset \Xi_m :=
\{|\zeta^m(0W)|_{\Tk}:\ 0W1 \in \Lk\}.
$$
We have
$$
|\zeta^m(0W)|_{\Tk} = \bt \ell(\zeta^m(0W)) = \bt M^m \ell(0W).
$$
Thus, 
$$
\Xi_m = \{\bt M^m z:\ z \in \Phi' \}
$$
by (\ref{Phipr}). Since $t_1/t_0$ is irrational and $M$ is an 
invertible integer
matrix, it is easy to see that $\bt M^m$ is a vector with
irrational ratio of components. (If $M$ were not invertible, then 
$\theta_2$ would be zero and $b(n)-a(n)$ would be bounded.) 
Let $\gamma_m$ denote this ratio.
Up to an overall scale, $\Xi_m$ is a projection of $\Phi'$
onto a line with slope $\gamma_m$. By Lemma~\ref{lem-elem2},
$\Phi' \supset \Z^2_+ \cap int(\Phti)$. By assumption and
Lemma~\ref{lem-elem22}, the width of $\Phti$
in the direction $(-\gamma_m,1)$ 
tends to infinity, which implies the desired statement. \qed

\begin{prop} \label{suff-2}
Let $\zeta$ be a primitive substitution on two letters.
If conditions (\ref{reprime}) and (\ref{excess}) are met,
then the $\Z$ action on $X_\zeta$ is topologically mixing. 
\end{prop} 

{\em Proof.}  Let $W_1, W_2$ be any allowed words. We need to show
that there exists $N\in \Nat$ such that for all $k\ge N$ there is an
allowed word $W_1 V W_2$ with $|W_1 V W_2|=k$. Since $\zeta$ is a primitive
substitution, there exists $m$ such that $\zeta^m(0)$ contains $W_1$
and $\zeta^m(1)$ contains $W_2$. Thus, it is enough to prove the claim
for $W_1 = \zeta^m(0)$ and $W_2 = \zeta^m(1)$.

By assumption (\ref{reprime}), we can find integers $r,s$ such that
\be 
\label{gcd} r |\zeta^m(0)| + s |\zeta^m(1)| = 1.  
\ee 
Thus, if two words $V$ and $V'$ have $\ell(V)-\ell(V')=
(r, s)^T$,
then $|\zeta^m(V)| - |\zeta^m(V')| = 1$.
We will demonstrate the existence of a sequence of words
$V_n$, each beginning with 0 and ending with 1, such that every
integer is within a bounded error of a length $|\zeta^m(V_n)|$. Once
the width of $\tilde \Phi$ in the $(-s/r, 1)$ direction is large enough,
there will exist other words $V_{n,j}$, also beginning with 0 and
ending with 1, whose population vectors differ from that of $V_n$ by
$(jr, js)^T$, 
so that $|\zeta^m(V_{n,j})| = |\zeta^m(V_n)| + j$.  For
every sufficiently large integer $k$, we can then pick $n$ and $j$
such that $|\zeta^m(V_{n,j})|=k$.

For each $n$ sufficently large that $b(n)-a(n) \ge 2$, let 
$$
p_n:= \lfloor(a(n) + b(n))/2 \rfloor,\ \ \ q_n:= n +1 - p_n,
$$
and let $V_n$ be a word, beginning with 0 and ending with 1, with
population vector $(p_n, q_n)$.  Let $u_n = |\zeta^m(V_n)|$.
It easily follows from Lemma~\ref{lem-elem1}(ii) that
$|p_{n+1}-p_n| \le 1$ and $|q_{n+1}-q_n| \le 1$, hence
$$
|u_{n+1} - u_n |\le |\zeta^m(0)| + |\zeta^m(1)|.
$$
Note also that $u_n \to \infty$. Thus, for any $k$ sufficiently large
we can find $n$ such that 
\be \label{star1}
u_n \le k \le u_{n+1} \le u_n + |\zeta^m(0)| + |\zeta^m(1)|.
\ee
If $k$ is sufficiently large, then $n$ is large as well, and we can make
sure that
\be \label{star2}
b(n) - a(n) \ge 5 + A,\ \ \ \mbox{where}\ \ A:=
\max\{|r|,|s|\} (|\zeta^m(0)| + |\zeta^m(1)|),
\ee
and hence that $b(n+1)-a(n+1) \ge 4 + A$. Now consider the points
$(p_n +r(k-u_n),q_n + s(k-u_n))$ and $(p_{n+1} + r(k-u_{n+1}), q_{n+1} +
s(k-u_{n+1}))$.  Since $k-u_n$ and $k-u_{n+1}$ have opposite signs,
one of these lies above (or on) the line traced by the points
$(p_n,q_n)$ and the other lies below (or on), and their distance is less than 
the width of $\tilde \Phi$.  Hence at least one of these points 
lies in the interior 
of $\tilde \Phi$ and there either exists
a word $V_{n,{k-u_n}}$ or a word $V_{n+1,k-u_{n+1}}$, beginning
with 0 and ending with 1.  Either way, $\zeta^m$ applied to this word has
length exactly $k$. 

\qed

This completes the proof of Theorem {\ref{thm2}}.


\section{Proof of Proposition \ref{prop1}}

Throughout this section we assume that $|\th_2|>1$.
Since the system only
depends on $u$, we can square the substitution and assume that $\theta_2 > 1$.
We begin with simple bounds on the numbers of 0s and 1s in substituted 
letters.  From these we show that $b(n)-a(n)$ is bounded above by a constant
times $n^\alpha$, where $\alpha = \log|\th_2|/\log\th_1$.
  These upper bounds, combined with an ``intermediate
value theorem'' argument, then yield lower bounds on $b(n)-a(n)$.

\begin{lemma} \label{lem-elem3}
There exist 
positive constants $L_1,L_1',L_2,L_2'$ such that for all $k\ge 1$,
\be \label{eq-new1}
L_1 \th_1^k \le \ell_0(\zeta^k(i)) \le L_1' \th_1^k,\ \ \ i=0,1,
\ee
and
\be \label{eq-new2}
L_2 |\th_2|^k \le
|\pigas(\ell(\zeta^k(i)))| \le L_2' |\th_2|^k = L_2'\th_1^{k\alpha},\ \ \ 
i=0,1.
\ee
Moreover, if $\th_2>0$, then $\pigas(\ell(\zeta^k(0)))<0$ and
$\pigas(\ell(\zeta^k(1)))>0$.
\end{lemma}

{\em Proof.} Fix $i\in \{0,1\}$ and
write $\ell(i) = a_1^{(i)} e_1 + a_2^{(i)} e_2$.  By (\ref{matr}),
$$
\ell(\zeta^k(i)) = M^k(\ell(i)) = a_1^{(i)}\theta_1^k e_1 +
a_2^{(i)}\theta_2^k e_2.
$$
The estimates (\ref{eq-new1}) hold since $|\theta_2| < \theta_1$
(note that $a_1^{(i)}>0$ since $\ell(i)$ is a non-negative vector).
The estimates (\ref{eq-new2}) follow from the fact that $\pigas(e_1)=0$.
To verify the last statement, we observe that $\ell(0) = (1,0)^T$, so 
$a_2^{(i)}<0$ by our choice of eigenvectors. On the other hand,
$\ell(1) = (0,1)^T$, so $a_2^{(i)} > 0$.
\qed

Thanks to Lemma \ref{lem-elem22}, the following estimate is 
tantamount to an upper bound on $b(n)-a(n)$.

\begin{lemma} \label{lem-accord}
Suppose that $|\theta_2|>1$. Then
there is a constant $C_2 > 0$, depending on $\zeta$, such that for any 
$w \in \Phi$ we have
$$
|\pigas(w)| \le C_2 w_x^\alpha.
$$
\end{lemma}

{\em Proof.}
A vector $w$ is in $\Phi$ if and only if there is a word $V\in \Lk$ 
such that $\ell(V) = w$. Since $V$ occurs in $u=\zeta(u)$,
we can write $V$ in the following form
(sometimes called the ``accordion form"):
$$
V = s_1 \zeta(s_2) \ldots \zeta^{k-1}(s_{k-1}) \zeta^{k}(s_k)
\zeta^k(p_k) \zeta^{k-1}(p_{k-1})\ldots\zeta(p_2)p_1
$$
where $s_j$ and $p_j$ are respectively suffixes and prefixes (possibly empty)
of the words $\zeta(i)$, $i \in \Ak$. Note that the number of possible words
$s_j$ and $p_j$ is finite (at most $|\zeta(0)|+|\zeta(1)|$).
We have for any word $W$,
$$
\ell(\zeta^n(W))  = M^n(\ell(W)) = b_1 \theta_1^n e_1 + 
b_2 \theta_2^n e_2\ \ \ \mbox{for all}\ \ n\ge 1,
$$
for some constants $b_1, b_2$. 
Since $\pigas(e_1) = 0$, we obtain
$$
|\pigas(w)| \le \const \cdot \sum_{j=0}^k |\th_2|^j < \const\cdot 
|\theta_2|^k/(|\th_2|-1),
$$
with the constant depending on the substitution, but not on $w$.
On the other hand, $w_x = \ell_0(V) \ge \const'\cdot \theta_1^k$
by (\ref{eq-new1}), and the
statement of the lemma follows, since $\theta_1^\alpha = |\theta_2|$.
\qed

\begin{prop}\label{prop3}
There is a constant $C_3 > 0$ such that for any $w\in \Z^2_+$
satisfying
\begin{equation} \label{est1}
-C_3 w_x^\alpha < \pigas(w) < 0,
\end{equation}
we have $w \in \Phi$.
\end{prop}

{\em Proof of Proposition \ref{prop1} assuming Proposition \ref{prop3}.}
The latter implies that the width of $\Phti$ in the direction of
$(-\gam,1)^T$ at the level $r$ is bounded below by $\const\cdot r^\alpha$
for any $r>0$. Then Lemma~\ref{lem-elem22}, together with the observation
preceding it, implies the desired estimate. \qed

\medskip

For the proof of Proposition \ref{prop3} we need
a simple geometric lemma, a kind of
``Intermediate value theorem."
We say that a point $z\in \Z^2_+$ is above 
(resp. below) $\Gamma$ if $z \not \in \Gamma$ and there are
no points of $\Gamma$ directly above (resp. below) it.
There is a natural linear order on $\Gamma$, with $(0,0)$ being the minimal
element, and the $(n+1)$st point always one unit above or to the right from
the $n$th point. Thus, every integer lattice point in the first quadrant either
belongs to $\Gamma$, or is above $\Gamma$, or is below $\Gamma$.

\begin{lemma} \label{lem-ivt}
If $w\in \Z^2_+$ and there exist $z,z' \in \Gamma$ such that $z+w$ is below
$\Gamma$ and $z'+w$ is above $\Gamma$, then $w \in \Gamma - \Gamma$.
\end{lemma}

The proof is straightforward, 
since when we move from a point $z$ on $\Gamma$ to the
next one, $z+w$ cannot ``jump"
from being below $\Gamma$ to being above $\Gamma$, or vice versa.

\medskip

{\em Proof of Proposition \ref{prop3}}.
Let $w$ be a vector satisfying (\ref{est1}). First we claim that there exists
$z \in \Gamma$ such that $z+w$ is on or below $\Gamma$. Suppose this is
not true. Since $z_0:=(0,0)\in \Gamma$,
we have that $w=z_0+w$ is above $\Gamma$, hence
there exists $z_1 \in \Gamma$ with
$$
(z_1)_x = w_x,\ \ \ \ (z_1)_y < w_y.
$$
Since $z_1+w$ is above $\Gamma$, there exists $z_2 \in \Gamma$ with
$(z_2)_x = 2 w_x$ and $(z_2)_y < 2 w_y$. 
Iterating this procedure, we obtain a sequence $z_m \in \Gamma$ such that
$$
(z_m)_x = mw_x,\ \ \ \ (z_m)_y < mw_y.
$$
This implies that
$$
\pigas(z_m) < m \pigas(w) = -\frac{|\pigas(w)|}{w_x} (z_m)_x,
$$
while Lemma \ref{lem-accord} implies that $\pigas(z_m) \ge -C_2(z_m)_x^\alpha$.
This is a contradiction for $m$ sufficiently large.
Notice that here we only used that $\pigas(w)<0$.

\medskip

In view of Lemma~\ref{lem-ivt}, it
remains to prove that there exists $z'\in \Gamma$ such that $z'+w$
is on or above $\Gamma$. Suppose there is no such $z'$.
For $k\ge 1$ let $\xi_k := \ell(\zeta^k(0))$, which is a point on $\Gamma$
by definition. Recall that in Lemma~\ref{lem-elem3} we showed that
$
L_1 \th_1^k \le (\xi_k)_x \le L_1' \th_1^k
$
and
\be \label{eq-new21}
\pigas(\xi_k) \le -L_2 \th_2^k = -L_2\th_1^{k\alpha}.
\ee
Let
\be \label{def-L4}
L_3:= \left(\frac{L_2}{2C_2\th_1^\alpha}\right)^{1/\alpha},
\ee
where $C_2$ is from Lemma~\ref{lem-accord}, and
\be \label{def-c1}
C_3:= C_2L_3/L_1'.
\ee
Further, let $k$ be the integer satisfying
\be \label{eq-new3}
L_3 \th_1^k \le w_x < L_3 \th_1^{k+1}.
\ee
We can assume that $w_x$ is sufficiently large, so that $k\ge 1$, since
for small $w_x$ the statement of the proposition is obviously true,
perhaps with a different constant.
Since $z_0:=(0,0)\in \Gamma$, the point $w=z_0+w$ is below $\Gam$
by our assumption, so there exists $z_1 \in \Gam$ such
that
$$
(z_1)_x = w_x,\ \ \ \ (z_1)_y > w_y.
$$
Iterating this, we obtain $z_m \in \Gam$ such that
\be \label{eq-new3.5}
(z_m)_x = mw_x,\ \ \ \ (z_m)_y > mw_y.
\ee
Let $m$ be the largest integer satisfying
$
mw_x \le (\xi_k)_x, 
$
so that 
\be \label{eq-new4}
(\xi_k - mw)_x < w_x.
\ee
We have
\be \label{eq-new5}
\pigas(\xi_k) = \pigas(\xi_k-z_m) + \pigas(z_m).
\ee 
Now,
\be \label{eq-new6}
|\pigas(\xi_k-z_m)| \le C_2 |(\xi_k- z_m)_x|^\alpha = C_2(\xi_k - mw)_x^\alpha
< C_2 w_x^\alpha < C_2 L_3^\alpha \th_1^{(k+1)\alpha},
\ee
by Lemma~\ref{lem-accord}, (\ref{eq-new4}) and (\ref{eq-new3}).
On the other hand,
\begin{eqnarray*}
\pigas(z_m) & > & m \pigas(w)\\ &  > & -C_3 m w_x^\alpha \\ &  > &
-C_3 m (L_3 \th_1^{k+1})^\alpha \\ & \ge & -C_3 \frac{(\xi_k)_x}{w_x}
L_3^\alpha \th_1^{(k+1)\alpha} \\ &  \ge & -C_3 \frac{L_1'\th_1^k}{L_3\th_1^k}
L_3^\alpha \th_1^{(k+1)\alpha} \\ & = &
-C_3 L_1' L_3^{\alpha-1} \th_1^{(k+1)\alpha} = -C_2 L_3^\alpha
\th_1^{(k+1)\alpha}.
\end{eqnarray*}
Above, we used (\ref{eq-new3.5}) in the 1st line, (\ref{est1}) in the
2nd line, (\ref{eq-new3}) in the 3d line,
the definition of $m$ in the 4th line,
(\ref{eq-new1}) and (\ref{eq-new3}) in the 5th line, and
(\ref{def-c1}) in the last line.
Combined with (\ref{eq-new5}) and (\ref{eq-new6}), the last inequality
yields
$$
\pigas(\xi_k) > -2C_2 L_3^\alpha \th_1^{(k+1)\alpha},
$$
contradicting (\ref{eq-new21}) in view of (\ref{def-L4}). \qed


\section{Proofs of other results, concluding remarks, and open questions}
\label{examplesection}

\begin{example}
Consider the ``period doubling'' substitution on two
letters $\{0,2\}$: $\zeta(0) = 02$, $\zeta(2) = 00$, whose
matrix $\left [ \begin{array}{ll} 1 & 2 \\ 1 & 0\end{array} \right ]$ 
has eigenvalues 2 and $-1$. 

Since each $\zeta(i)$ begins with 0 and has length 2, every word of
length $2n$ is either $\zeta$ applied to a word of length $n$, or is obtained
from such a word by deleting the initial 0 and adding a 0 at the end. In
particular, the population vector of each word of length $2n$ is $M$ times 
the population vector of a word of length $n$, from which we infer that
$b(2n)-a(2n) = b(n)-a(n)$. In particular, $b(2^n)-a(2^n)=b(1)-a(1) = 1$. 
A word of length $2^n$ has either $\lfloor 2^{n+1}/3 \rfloor$ 0's and
$\lfloor 2^n/3 \rfloor + 1$ 2's or $\lfloor 2^{n+1}/3 \rfloor + 1$ 0's and
$\lfloor 2^n/3 \rfloor$ 2's.

This subshift is neither weak mixing nor topologically mixing, as
condition (\ref{reprime}) is not met.  Since 
$\liminf (b(n)-a(n)) = 1$, the $\R$ action on a tiling space based on this
substitution  will not be topologically mixing.  However,
the $\R$ action will be weak mixing 
if the ratio of tile lengths is irrational \cite{CS}.
\end{example}
 
\begin{example} ( = Example 2.2). Consider the ``modified period doubling
substitution'' (MPD)
$\zeta(0)= 011,\ \zeta(1)=0$. 
The matrix of the 
substitution is $M= \left[ \begin{array}{ll} 1 & 1 \\ 2 & 0 
\end{array} \right]$, with eigenvalues $\theta_1 = 2$, $\theta_2 = -1$.
Condition (\ref{reprime}) is easily seen to hold, since the only prime that 
divides the determinant of $M$ is 2, and $(1,1) M^n = (1,1) \pmod{2}$. 

This subshift
is obtained from the period doubling subshift by replacing each
2 with a pair of 1's. The tiling spaces built from this subshift are 
topologically conjugate to the tiling spaces of the period doubling subshift,
if we take the same values of $t_0$ and let $t_1 = t_2/2$. In particular,
the $\R$ action is weak mixing if $t_0/t_1$ is irrational, but is not
topologically mixing. This implies that $\liminf (b(n)-a(n))$ is finite,
so the $\Z$ action is not topologically mixing. 

We can obtain precise bounds on $\liminf (b(n)-a(n))$ from the correspondence
with the period doubling subshift.  Every MPD word is obtained from a
period-doubling word by replacing each 2 with 11 and then possibly truncating
leading or trailing 11's to a single 1.  From period-doubling words of length
$2^n$ we obtain MPD words of length $\lfloor 2^{n+2}/3 \rfloor$ with
$\lfloor 2^{n+1}/3 \rfloor +1$ 0's, of length  $\lfloor 2^{n+2}/3 \rfloor$ with
$\lfloor 2^{n+1}/3 \rfloor$ 0's, or of length  $\lfloor 2^{n+2}/3 \rfloor+1$ 
with $\lfloor 2^{n+1}/3 \rfloor$ 0's. One can also obtain MPD words of
length  $\lfloor 2^{n+2}/3 \rfloor$ with $\lfloor 2^{n+1}/3 \rfloor -1$ 0's 
from period-doubling words of length $2^n-1$.  
Thus $b(\lfloor 2^{n+2}/3 \rfloor) = \lfloor 2^{n+1}/3 \rfloor +1$ 
and $a(\lfloor 2^{n+2}/3 \rfloor) = \lfloor 2^{n+1}/3 \rfloor -1$, 
so $\liminf b(n)-a(n)$ is at most 2. 
\end{example}

Fix $t_0$ and $t_1$ with irrational ratio, and consider any tiling of
$\R$ associated with the MPD substitution. Denote by $\La$ the
set the endpoints of its tiles. It provides an interesting example of
a uniformly discrete, relatively dense subset of the real line.
The arithmetic difference $\La-\La$ is just the projection of $\Ga-\Ga$
along an irrational direction. Since the width of $\widetilde{\Phi}$
has infinite $\limsup$ and finite $\liminf$, we obtain that $\La-\La$
is not uniformly discrete, but it is not eventually dense in $\R$ either.
A subset of $\R$ is said to be {\em Meyer} if it is relatively dense and
its arithmetic self-difference
is uniformly discrete. Thus, we obtain the following. 

\begin{cor} \label{cor-meyer} There exists a ``binary'' (i.e.\ with two
possible distances between consecutive points) 
non-Meyer set $\Lambda \subset \R$
such that $\Lambda-\Lambda$ is not eventually dense.
\end{cor}

\noindent
{\bf Open questions and directions for further research.}

\smallskip

{\bf 1.} Is there a combinatorial description of
substitutions on two symbols that have
$|\theta_2|=1$ and satisfy (\ref{excess})? Are there primitive 
integer  matrices
with the second eigenvalue of magnitude one, which determine
topological mixing (or lack thereof) for every substitution having this
matrix?

\smallskip

{\bf 2.} What happens for substitutions on more than two symbols? Of course,
there are trivial examples arising from the fact that a substitution may be
recoded using a new, larger alphabet. With some caution, we can put
forward the following {\em Conjecture: for a primitive substitution
on any number of symbols, assuming none of the eigenvalues has magnitude
equal to one, topological mixing is equivalent to weak mixing.}
This concerns both the $\Z$ and $\R$ actions associated with the
substitution.

\smallskip

{\bf 3.} What happens in higher dimensions? We do not know any examples of
topologically mixing substitution $\Z^d$ or $\R^d$ actions for $d\ge 2$.
In particular, is the pinwheel tiling system studied
by Radin and coauthors (see \cite{Radin}) topologically mixing?
A well-known open problem is to determine whether the pinwheel system is
strongly mixing. Checking topological mixing may be a more realistic goal.

\medskip

\noindent {\bf Acknowledgments.} Thanks to Natalie Priebe Frank for helpful
comments.
This work was largely inspired by a question
of Robbie Robinson at the BIRS Workshop on Joint Dynamics in Banff, Canada,
in June 2003. L.S. and B.S. would like to thank Michael Baake,
the organizer of the August 2003 Workshop on Mathematics of Aperiodic Order in
Greifswald, Germany, where part of this work was done.
This work was completed while R.K. was visiting Princeton University.


\bibliographystyle{amsplain}

\end{document}